\documentclass[reqno]{amsart}

\usepackage{color}

\theoremstyle{plain}
\newtheorem{theorem}{Theorem} [section]
\newtheorem{corollary}[theorem]{Corollary}
\newtheorem{lemma}[theorem]{Lemma}
\newtheorem{proposition}[theorem]{Proposition}

\theoremstyle{definition}
\newtheorem{definition}[theorem]{Definition}

\newtheorem{standing}[theorem]{Standing Assumptions}
\newtheorem{remark}[theorem]{Remark}

\newtheorem{example}[theorem]{Examples}

\def\C{{\mathbb{C}}}
\def\N{{\mathbb{N}}}
\def\R{{\mathbb{R}}}
\def\Z{{\mathbb{Z}}}

\def\A{{\mathcal{A}}}

\def\clspan{{\overline{\mathrm{span}}}}

\newcommand{\p}{\varphi}
\newcommand{\w}{\omega}
\newcommand{\e}{\varepsilon}
\newcommand{\f}{\alpha}

\newcommand{\g}{\gamma}
\newcommand{\dd}{\delta}
\newcommand{\G}{\Gamma}
\newcommand{\W}{\Omega}
\newcommand{\Dd}{\Delta}
\newcommand{\T}{\mathcal{T}}

\def\Span{{\text{\rm span}}}  
\def\subset{\subseteq}

\begin{document}

\title{Shift Invariant Spaces on LCA Groups}
\author[C. Cabrelli and V. Paternostro]{Carlos Cabrelli and Victoria Paternostro}

\address{\textrm{(V.Paternostro)}
Departamento de Matem\'atica,
Facultad de Ciencias Exac\-tas y Naturales,
Universidad de Buenos Aires, Ciudad Universitaria, Pabell\'on I,
1428 Buenos Aires, Argentina and
CONICET, Consejo Nacional de Investigaciones
Cient\'ificas y T\'ecnicas, Argentina}
\email{vpater@dm.uba.ar}

\address{\textrm{(C. Cabrelli)}
Departamento de Matem\'atica,
Facultad de Ciencias Exac\-tas y Naturales,
Universidad de Buenos Aires, Ciudad Universitaria, Pabell\'on I,
1428 Buenos Aires, Argentina and
CONICET, Consejo Nacional de Investigaciones
Cient\'ificas y T\'ecnicas, Argentina}
\email{cabrelli@dm.uba.ar}

 \thanks{The research of
 C. Cabrelli and V. Paternostro is partially supported by
Grants: ANPCyT, PICT 2006--177, CONICET, PIP 5650, UBACyT X058 and X108.
}

\subjclass[2000]{Primary 43A77; Secondary 43A15.}

\keywords{Shift-invariant spaces, translation invariant spaces,
LCA groups, range functions, fibers}

\date{\today}
\maketitle

\begin{abstract}
In this article we extend the theory of shift-invariant spaces to the context of LCA groups.
We introduce the notion of $H$-invariant space for a countable discrete subgroup $H$ of an
LCA group $G$,
and show that the concept of range function and the techniques of fiberization are valid in this context. As a consequence of this generalization we prove characterizations of frames and
Riesz bases of these spaces extending previous results, that were known for $\R^d$ and the lattice $\Z^d.$
\end{abstract}

\section{Introduction}
A {\it shift-invariant space} (SIS) is a closed subspace of $L^2(\R)$ that is invariant under translations by integers.
The Fourier transform of a shift-invariant space is a closed subspace that  is invariant under integer modulations (multiplications by complex exponentials of integer frequency). Spaces that are invariant under integer modulations are called  {\it doubly invariant spaces}.
Every result on doubly invariant spaces can be translated to an equivalent result in shift-invariant spaces via the Fourier transform.
Doubly invariant spaces have been studied in the sixties by Helson \cite{helson} and also by Srinivasan \cite{S1}, \cite{S2}, in the context of  operators related to harmonic analysis.

Shift-invariant spaces are very important in applications and the theory had a great development in the
last twenty years, mainly in approximation theory, sampling, wavelets, and frames. In particular they serve as
models in many problems in signal and image processing.

In order to understand the  structure of doubly invariant spaces, Helson introduced  the notion of {\it range function}.
This  became an essential tool in the modern development of the theory. See \cite{dBDR94a}, \cite{dBDR94b}, \cite{RS95} and \cite{B}.

Range functions characterize completely shift-invariant spaces and provide a series of techniques
known in the literature as {\it fiberization} that allow to have a different view and a deeper insight of these spaces.

Fiberization techniques are very important in the class of {\it finitely generated} shift-invariant spaces.
A key feature of these spaces is that they can be generated by the integer translations of a finite number of functions. Using range functions  allows us to translate  problems on finitely
generated shift-invariant spaces,  into problems of linear algebra (i.e. finite dimensional problems).

Shift-invariant spaces generalize very well to several variables where the invariance is understood to be under the group $\Z^d$.

When looking carefully at the theory it becomes apparent that it is strongly based on the additive group operation of $\R^d$ and the  action of the subgroup $\Z^d.$

It is therefore interesting to see if the theory can be set in a context of  general {\it locally compact
abelian} groups (LCA groups).
The locally compact abelian group framework has several advantages.
First because it  is important to have a valid theory  for the classical groups such as $ \Z^d, \mathbb{T}^d$ and
$ \Z_n$. This will be crucial particularly in applications, as in the case of the generalization of the Fourier Transform to LCA groups and also Kluvanek's theorem, where the Classical Sampling theorem is extended to this general  context, (see \cite{Klu65}, \cite{Dodson}).

On the other side,  the LCA groups setting,  unifies a number of different results into a general framework with a concise and elegant notation. This fact enables us to visualize hidden relationships  between the different components of the theory, what, as a consequence, will translate in a deeper and better understanding of shift-invariant spaces, even in the case of the real line.

In this paper we develop the theory of shift invariant spaces in LCA groups.
Our emphasis will be  on range functions and fiberization techniques. The order of the subjects  follows mainly the
treatment  of Bownik in $\R^d$, \cite{B}.
In \cite{KR} the authors study, in the context of LCA groups, {\it principal} shift-invariant spaces,
that is, shift-invariant spaces generated by one single function. However they don't
develop the general theory.

This article is organized in the following way.
In Section \ref{sec-1} we give the necessary  background on LCA groups and set the basic notation. In Section \ref{sec-2} we state our standing assumptions and prove the characterizations of $H$-invariant spaces using range functions. We apply these results in Section \ref{sec-3}
to obtain a characterization of frames and Riesz bases of $H$-translations.

\section{Background on LCA Groups}\label{sec-1}

In this section we review some basic known results from the theory of LCA groups, that
we  need for the remainder of the article.  In this way we set the notation that we will
use in the following sections. Most proofs of the results are omitted unless it is considered necessary.
For details and proofs  see \cite{rudin},\cite{hewitt-ross-1}, \cite{hewitt-ross-2}.

\subsection{LCA Groups}\noindent

In this section we review some basic known results from the theory of LCA groups, that
we  need for the remainder of the article.  In this way we set the notation that we will
use in the following sections. Most proofs of the results are omitted unless it is considered necessary.
For details and proofs  see \cite{rudin},\cite{hewitt-ross-1}, \cite{hewitt-ross-2}.

\subsection{LCA Groups}\noindent

Throughout this article, $G$ will denote a locally compact abelian, Hausdorff group (LCA) and
 $\G$ (or $\widehat{G}$)  its dual group. That is,
$$\Gamma=\big\{\gamma:G\rightarrow \C: \gamma \,\,\textrm{is a continuous character of }\, G\big\},$$
where a character is a function such that:
\renewcommand{\theenumi}{\alph{enumi}}
\begin{enumerate}
\item[(a)] $|\g(x)|=1,\,\,\forall \,\,x\in G$.
\item[(b)] $\g(x+y)=\g(x)\g(y),\,\,\forall\,\, x,y \in G$.
\end{enumerate}
Thus,  characters  generalize the exponential functions
$\g_t(y)=e^{2\pi i ty}$, from the case  $G=(\R,+)$.

Since in this context, both the algebraic and topological structures coexist, we will say that
two groups $G$ and $G'$ are {\it topologically isomorphic} and we will write $G \approx  G'$,
if there exists a topological isomorphism
from $G$ onto $G'$. That is, an algebraic isomorphism which is a homeomorphism as well.

The following theorem states some important facts about LCA groups. Its proof can be found in \cite{rudin}.

\begin{theorem}
Let $G$ be an LCA group and  $\G$ its dual. Then,
\begin{enumerate}
\item [(a)]The dual group $\G$, with the operation $(\g+\g')(x)=\g(x)\g'(x)$, is an LCA group.
\item [(b)]The dual group of $\G$ is topologically isomorphic to $G$, with the identification
$x\in G\leftrightarrow \phi_x\in\widehat{\G},$ where  $\phi_x(\g):=\g(x)$.
\item [(c)]$G$ is discrete (compact) if al only if $\G$ is compact (discrete).
\end{enumerate}
\label{teo-dualidad}
\end{theorem}

As a consequence of item $(b)$ of Theorem \ref{teo-dualidad}, it is convenient to use the notation $(x,\g)$ for the complex number $\g(x)$, representing the character $\g$ applied to $x$ or the character $x$ applied to $\g$.

Next we list the most basic examples that are  relevant to Fourier analysis.
As usual, we identify the interval $[0,1)$ with the torus  $\mathbb{T}=\{z\in\C\,:\,|z|=1\}.$
\begin{example}
\renewcommand{\theenumi}{\Roman{enumi}}
\noindent
\begin{enumerate}
\item[(I)] In case that $G=(\R^d,+)$, the dual group $\G$   is also $(\R^d,+)$, with the identification
$x\in \R^d\leftrightarrow\g_x\in\G$, where $\g_x(y)=e^{2\pi i \langle x,y\rangle}$.
\item[(II)] In case that $G=\mathbb{T}$,
its dual group is topologically isomorphic to $\Z$, identifying
each $k\in\Z$ with $\g_k\in\G$, being $\g_k(\w)=e^{2\pi i k\omega}$.
\item[(III)] Let  $G=\Z$. If $\g\in\G$, then $(1,\g)= e^{2\pi i\f}$ for same  $\f\in\R$.
Therefore, $(k,\g)= e^{2 \pi i\f k}$. Thus, the complex number $ e^{2 \pi i\f}$ identifies the  character $\g.$
This proves that $\G$ is $\mathbb{T}$.
\item[(IV)] Finally, in case that $G=\Z_n$, the dual group  is also $\Z_n$.
\end{enumerate}
\end{example}
\renewcommand{\theenumi}{\roman{enumi}}

Let us now consider  $H\subseteq G$,  a closed subgroup of an LCA group $G$. Then, the quotient $G/H$ is a regular (T3)
topological group.
Moreover,  with the quotient topology, $G/H$ is  an LCA group and
if $G$ is second countable,  the quotient $G/H$  is also second countable.

For an LCA group $G$  and   $H\subset G$  a  subgroup of $G$, we define the subgroup $\Delta$
of $\G$ as follows:
$$\Delta=\big\{\gamma \in \Gamma : (h,\g)=1, \,\,\forall\,\, h\in H\big\}.$$
This subgroup is called the \textit{annihilator of $H$}.
Since each character in $\G$ is a continuous function on $G$,
$\Dd$ is a closed subgroup of $\Gamma$. Moreover, if $H\subset G$ is a closed subgroup  and $\Dd$ is the annihilator of $H$, then $H$ is the annihilator of $\Dd$ (see \cite[Lemma 2.1.3.]{rudin}).

The next result establishes duality relationships among the groups
$H$, $\Delta$, $G/H$ and $\Gamma/\Delta$.

\begin{theorem}
If $G$ is an LCA group and  $H\subset G$ is a closed subgroup of $G$, then:
\begin{enumerate}
\item[(i)]\label{T1}  $\Delta$ is topologically isomorphic to the dual group of $G/H$, i.e:
$\Dd\approx\widehat{(G/H)}$.
\item[(ii)]\label{T2}  $\Gamma/\Delta$ is topologically isomorphic to the dual group of $H$, i.e:
$\G/\Dd\approx \widehat{H}$.
\end{enumerate}
\label{teo-de-isos}
\end{theorem}

\begin{remark}
According to Theorem \ref{teo-dualidad}, each element of $G$ induces one character in $\widehat{\G}$.
In particular, if $H$ is a closed
subgroup of $G$, each $h\in H$ induces a character that has the additional property of being $\Dd$-periodic.
That is, for every $\dd\in\Dd$, $(h,\g+\dd)=(h,\g)$ for all $\g\in\G$.
\label{caracteres-delta-periodicos}
\end{remark}

The following definition will be useful throughout this paper. It agrees with the
one given in \cite{ceros-zak}.
\begin{definition}
Given $G$ an LCA group, a  {\it uniform lattice $H$ in $G$} is a discrete subgroup of $G$ such that the quotient group $G/H$
is  compact.
\end{definition}

The next theorem points out a number of relationships which occur among $G$, $H$, $\Gamma$, $\Delta$ and
their respective quotients.

\begin{theorem}\label{theo-standing}
\renewcommand{\theenumi}{\arabic{enumi}}
Let
$G$ be a second countable LCA group. If
$H\subset G$ is a countable (finite or countably infinite) uniform lattice,
the following properties hold.

\begin{enumerate}
\item[(1)] $G$ is separable
\item[(2)] $H\subset G$ is  closed.
\item[(3)] $G/H$ is second countable  and metrizable.
\item[(4)] $\Dd\subset\G$, the annihilator of $H$, is closed, discrete and countable.
\item[(5)] $\widehat{H}\approx\G/\Dd$ and $\widehat{(G/H)}\approx\Dd$.
\item[(6)] $\G/\Dd$ is a compact group.
\end{enumerate}
\end{theorem}

Note that in particular, this theorem states that $\Dd$ is a countable uniform lattice
in $\G$.

\renewcommand{\theenumi}{\roman{enumi}}
\subsection{Haar Measure on LCA groups}\noindent

On every LCA group $G$, there exists a {\it Haar measure}. That is, a non-negative,
regular borel measure  $m_G$,
which is not identically zero and {\it trans-lation-invariant}.
This last property means that,
$$m_G(E+x)=m_G(E)$$
for every element $x\in G$ and every Borel set $E\subset G$. This measure is
unique up to constants, in the following sense: if $m_G$ and $m'_G$ are two Haar measures on $G$, then
there exists a positive constant $\lambda$  such that $m_G=\lambda \,m'_G$.

Given a Haar measure $ m_G$ on an LCA group $G$,
the integral over $G$ is translation-invariant in the sense that,
$$\int_G f(x+y)\,dm_G(x)=\int_G f(x)\,dm_G(x)$$
for each element $y\in G$ and for each Borel-measurable function $f$ on $G$.

As in the case of the Lebesgue measure, we can define the spaces $L^p(G,m_G)$,
that we will  denote as $L^p(G)$, in the following way
$$ L^p(G) = \big\{ f:G\rightarrow \C: f \text{ is measurable and  } \int_G |f(x)|^p\,dm_G(x)<\infty\big\}.$$
If $G$ is a second countable LCA group,  $L^p(G)$ is
separable, for all $1\leq p<\infty$.
We will focus here on the cases $p=1$ and $p=2$

The next theorem is a generalization of the periodization argument usually
applied in case $G=\R$ and $H=\Z$  (for details see \cite[Theorem 28.54]{hewitt-ross-2}).

\begin{theorem}
Let $G$ be an LCA group, $H\subset G$  a closed subgroup and  $f\in L^1(G)$. Then,
the Haar measures $m_G$, $m_H$ and $m_{G/H}$  can be chosen such that
$$\int_G f(x)\,dm_G(x)=\int_{G/H} \int_H f(x+h)\,dm_H(h)\,dm_{G/H}([x]),$$
where $[x]$ denotes the coset of $x$ in the quotient $G/H$.
\label{teo-coset-descomposition}
\end{theorem}

If $G$ is a countable discrete group, the integral of
$f\in L^1(G)$ over $G$, is determined by the formula
$$\int_G f(x)\,dm_G(x)=m_G(\{0\})\sum_{x\in G}f(x),$$
since, due to the translations invariance, $m_G(\{x\})=m_G(\{0\})$, for each element $x\in G$.

\begin{definition}
A {\it section} of $G/H$ is a set of representatives of this quotient.
That is, a subset $C$ of $G$ containing exactly  one  element of each coset. Thus,
each element $x\in G$ has a unique expression of the form
$x=c+h$ with $c\in C$ and $h\in H$.
\end{definition}

We will need later in the paper to work with  Borel sections.
The existence of Borel sections  is provided by the following
lemma (see \cite{ceros-zak} and \cite{seccion-boreliana}).

\begin{lemma}
Let $G$ be an LCA group and  $H$ a uniform lattice in $G$.
Then, there exists a section of the quotient $G/H$, which is Borel
measurable.

Moreover, there exists a section of $G/H$ which is relatively compact.
\label{seccion-medible}
\end{lemma}

A  section $C\subset G$ of $G/H$  is in  one to one correspondence with $G/H$ by the {\it cross-section} map
$\tau:G/H\rightarrow  C$, $[x]\mapsto [x] \cap C$.
Therefore, we can carry over the topological and algebraic structure of  $G/H$ to $C$. Moreover, if $C$ is a Borel section, $\tau:G/H\rightarrow  C$ is measurable with respect to the Borel $\sigma$-algebra in $G/H$ and the Borel $\sigma$-algebra  in $G$ (see \cite[Theorem 1 ]{seccion-boreliana}). Therefore, the set value function defined by $m(E)=m_{G/H}(\tau^{-1}(E))$ is well defined on Borel subsets of $C$.
In the next lemma, we will prove that this measure $m$ is equal to $m_G$ up to a constant.

\begin{lemma}\label{lema-medidas}
Let $G$ be an LCA group, $H$ a countable uniform lattice in $G$ and  $C$  a Borel section of $G/H$. Then,
for every Borel set $E\subset C$
$$m_G(E)=m_H(\{0\})m_{G/H}(\tau^{-1}(E)),$$
where $\tau$ is the cross-section map.

In particular,
$ m_G(C)=m_H(\{0\})m_{G/H}(G/H).$
\end{lemma}

\begin{proof}
According to Lemma \ref{seccion-medible}, there exists a relatively compact  section of $G/H$. Let us call it  $C'$. Therefore, if $C$ is any other Borel section of $G/H$, it must satisfy $m_G(C)=m_G(C')$.
Since $C'$ has finite $m_G$ measure, $C$ must have finite measure as well.

Now, take $E\subset C$ a Borel set. Using Theorem \ref{teo-coset-descomposition},
\begin{eqnarray*}
m_G(E)=\int _G \chi_E(x)\,dm_G(x)
&=&\int_{G/H}\int_H \chi_E(x+h)dm_H(h)\,dm_{G/H}([x])\\
&=&m_H(\{0\})\int_{G/H}\sum_{h\in H}\chi_E(x+h)\,dm_{G/H}([x])\\
&=&m_H(\{0\})\int_{G/H}\chi_{\tau^{-1}(E)}([x])\,dm_{G/H}([x])\\
&=&m_H(\{0\})m_{G/H}(\tau^{-1}(E)).
\end{eqnarray*}
\end{proof}

\begin{remark}\label{estructura-seccion}
Notice that $C$, together with the LCA group structure inherited by $G/H$ through $\tau$, has the Haar measure $m$.
We  proved that $m_G\vert_C$, the restriction of $m_G$ to $C$, is a multiple of $m$. It follows that
$m_G\vert_C$ is also a Haar measure on $C$.

In this paper we will consider  $C$ as an LCA group with
the structure inherited by $G/H$ and with the Haar measure $m_G$.
\end{remark}

A {\it trigonometric polynomial}  in an LCA group $G$ is a function of the form
$P(x)=\sum_{j=1}^na_j(x,\g_j),$
where $\g_j\in\G$ and $a_j\in\C$ for all $1\leq j \leq n$.

As a consequence of  Stone-Weierstrass Theorem, the following result holds, (see \cite{rudin}, page 24).

\begin{lemma}
If $G$ is a compact LCA group, then the trigonometric polynomials are dense in
$\mathcal{C}(G)$,
where $\mathcal{C}(G)$ is the set of all continuous complex-valued functions on $G$.
\label{poly-triginometricos-densos}
\end{lemma}

Another important property of characters in compact groups is the following. For its proof see proof of \cite[Theorem 1.2.5]{rudin}.

\begin{lemma}
Let $G$ be a compact LCA group and  $\G$ its dual.
Then, the characters of  $G$ verify the following orthogonality relationship:
$$\int_G(x,\g)\overline{(x,\g')}\,dm_G(x)=m_G(G)\dd_{\g\g'},$$
for all $\g,\g'\in\G$, where $\dd_{\g \g'}=1$ if $\g=\g'$  and $\dd_{\g \g'}=0$ if $\g\neq\g'$.
\label{compacto-ortogonal}
\end{lemma}

Let us now suppose that $H$ is a uniform lattice in $G$.
If $\G$ is the dual group of $G$
and $\Dd$ is the annihilator of $H,$ the following characterization of the characters of the  group $\G/\Dd$ will be useful to understand what follows.

For each $h\in H$, the function
$\g\mapsto (h,\g)$
is constant on the cosets $[\g]=\g+\Dd$. Therefore, it defines a character on $\G/\Dd.$
Moreover, each character on $\G/\Dd$ is of this form. Thus, this correspondence between $H$ and the characters of $\G/\Dd$, which is actually a topological isomorphism, shows the dual relationship established  in Theorem \ref{teo-de-isos}.

Furthermore, since $\G/\Dd$ is compact, we can apply Lemma \ref{compacto-ortogonal} to $\G/\Dd.$
Then, for $h\in H$, we have
\begin{equation}
\int_{\G/\Dd}(h,[\g])\,dm_{\G/\Dd}([\g])=\left\{\begin{array}{cc}
m_{\G/\Dd}(\G/\Dd) & \mbox{ if }  h=0 \\
0 & \mbox{ if }  h\neq 0
\end{array}.
\right.
\label{caracteres-ortogonales}
\end{equation}

\subsection{The Fourier transform on LCA groups}\noindent

Given a function  $f\in L^1(G)$ we define the  {\it Fourier Transform} of $f$, as
\begin{equation}
\widehat{f}(\g)=\int_G f(x)(x,-\g)\,dm_G(x),\quad\g\in\G.
\label{l1-fourier}
\end{equation}

\begin{theorem}
The Fourier transform  is a linear operator
from $L^1(G)$ into $\mathcal{C}_0(\G)$,
where $\mathcal{C}_0(\G)$ is the subspace of $ \mathcal{C}(\G)$ of functions  vanishing at infinite, that is, $f \in \mathcal{C}_0(\G)$ if  $f \in \mathcal{C}(\G)$ and for all $\e>0$
there exists a compact set $K\subset G$ with $|f(x)|<\e$ if $x\in K^c$.

Furthermore, $\land:L^1(G)\rightarrow\mathcal{C}_0(\G)$ satisfies
\begin{equation}\label{fourier-mono}
\widehat{f}(\g)=0\,\,\forall\,\g\in\G\,\Rightarrow\,\,f(x)=0\,\, a.e.\,\,x\in G.
\end{equation}
\end{theorem}

The Haar measure of the dual group $\G$ of $G$, can be normalized so that, for a specific class
of functions,
the following inversion formula holds (see \cite[Section 1.5]{rudin}),
$$f(x)=\int_{\G}\widehat{f}(\g)(x,\g)\,dm_{\G}(\g).$$

In the case that the Haar measures $m_G$ and $m_{\G}$ are normalized such that the inversion formula holds,
the Fourier transform on $L^1(G)\cap L^2(G)$ can be extended
to a unitary operator from $L^2(G)$ onto $L^2(\G)$, the so-called
Plancharel transformation. We also denote this transformation by "$\land$".

Thus, the Parseval formula holds
$$\langle f,g\rangle=\int_G f(x)\overline{g(x)}\,dm_G(x)=
\int_{\G}\widehat{f}(\g)\overline{\widehat{g}(\g)}\,dm_{\G}(\g)=\langle \widehat{f},\widehat{g}\rangle,$$
where $f,g\in L^2(G) $.

Let us now suppose  that $G$ is compact. Then $\G$ is discrete.
Fix  $m_G$ and  $m_{\G}$ in order that
the inversion formula
holds. Then, using the Fourier transform, we obtain that
\begin{equation}\label{ecu-compacto-discreto}
1=m_{\G}(\{0\})m_{G}(G).
\end{equation}

The following lemma is a straightforward consequence of Lemma \ref{compacto-ortogonal}, equation (\ref{caracteres-ortogonales}) and statement (\ref{fourier-mono}).

\begin{lemma}\label{lema-caracteres-base}
If $G$ is a compact  LCA group and its dual $\G$ is countable, then the characters $\{\g\,:\,\g\in\G\}$
form an orthogonal basis for $L^2(G)$.
\end{lemma}

For an  LCA group $G$ and a countable uniform lattice $H$ in $G$,
we will denote by $\W$ a Borel section of $\Gamma/\Dd$.  In the remainder of this paper we will identify  $L^2(\W)$
with the set $\{\p\in L^2(\G): \p=0\,\, a.e.\,\,\, \G\setminus\W\}$ and $L^1(\W)$ with the set
$\{\phi\in L^1(\G):\, \phi=0\,\,\, a.e.\,\,\, \G\setminus\W\}$.

Let us now define  the functions $\eta_h:\G\mapsto\C$,
as $\eta_h(\g)=(h,-\g)\chi_{\W}(\g)$.
Using Lemma \ref{lema-caracteres-base} we have:
\begin{proposition}
Let $G$ be an LCA group and $H$ a countable uniform lattice in $G$.
Then, $\{\eta_h\}_{h\in H}$ is an orthogonal basis for $L^2(\W)$.
\label{lema-bon}
\end{proposition}

\begin{remark}\label{obs-iso}
We can associate to each $\p\in L^2(\W)$,   a  function $\p'$
defined on  $\G/\Dd$ as
$\p'([\g])=\sum_{\dd\in\Dd}\p(\g+\dd)$.
The correspondence $\p\mapsto\p'$, is an isometric isomorphism up to a constant between
$L^2(\W)$ and $L^2(\G/\Dd)$, since
$$\|\p\|^2_{L^2(\W)}=m_{\Dd}(\{0\})\|\p'\|^2_{L^2(\G/\Dd)}.$$
\end{remark}

Combining the above remark, Proposition \ref{lema-bon}, and the relationships established
in Theorem \ref{teo-de-isos}, we obtain the following proposition, which will be very important
on the remainder of the paper.

\begin{proposition}
Let $G$ be an LCA group, $H$ countable uniform lattice
on  $G$, $\G=\widehat{G}$ and  $\Dd$ the annihilator of $H$.
Fix $\,\W$ a Borel section of $\,\G/\Dd$ and  choose
$m_H$ and $m_{\G/\Dd}$ such that the inversion formula holds.
Then
$$\|a\|_{\ell^2(H)}=
\frac{m_H(\{0\})^{1/2}}{m_{\G}(\W)^{1/2}}\|\sum_{h\in H}a_h\eta_h\|_{L^2(\W)},$$
for each  $a=\{a_h\}_{h\in H}\in\ell^2(H)$.
\label{lema-parseval}
\end{proposition}

\begin{proof}
Let $a\in\ell^2(H)$. So,
\begin{equation}
\|a\|_{\ell^2(H)}=\|\widehat{a}\|_{L^2(\G/\Dd)},
\label{ecu-1-parseval}
\end{equation}
since $\widehat{H}\approx\G/\Dd$ and therefore $\land:H\rightarrow \G/\Dd$.

Take $\p(\g)=\sum_{h\in H}a_h(h,-\g)\chi_{\W}(\g)$. Then, by Proposition \ref{lema-bon}, $\p\in L^2(\W)$.
Furthermore,  $\p'([\g])=\p(\g)$, a.e. $\g\in\W$. So, as a consequence
of Remark \ref{obs-iso}, we have
\begin{equation}
\|\p'\|^2_{L^2(\G/\Dd)}=\frac1{m_{\Dd}(\{0\})}\|\p\|^2_{L^2(\W)}.
\label{ecu-2-parseval}
\end{equation}

Now,
$\widehat{a}([\g])=m_H(\{0\})\sum_{h\in H}a_h(h,-[\g]).$
Therefore, substituting in equations (\ref{ecu-1-parseval}) and (\ref{ecu-2-parseval}),
$$\|a\|_{\ell^2(H)}=\frac{m_H(\{0\})}{m_{\Dd}(\{0\})^{1/2}}\|\p\|_{L^2(\W)}.$$

Finally, since $m_{\G}(\W)=m_{\Dd}(\{0\})\,m_{\G/\Dd}(\G/\Dd)$, using (\ref{ecu-compacto-discreto})
we have that
$$\frac{m_H(\{0\})}{m_{\Dd}(\{0\})^{1/2}}=\frac{m_H(\{0\})^{1/2}}{m_{\G}(\W)^{1/2}},$$
which completes the proof.
\end{proof}

We finish this section with a result which is a consequence of statement (\ref{fourier-mono}) and Theorem \ref{teo-coset-descomposition}.

\begin{proposition}
Let $G$, $H$ and $\W$ as in Proposition \ref{lema-parseval}.
If $\phi\in L^1(\W)$
and $\widehat{\phi}(h)=0$ for all $h\in H$, then $\phi(\w)=0$ a.e. $\w\in\W$.
\label{prop-fourier-H-0}
\end{proposition}

\section{$H$-Invariant Spaces}\label{sec-2}
In this section we extend the theory of shift-invariant
spaces in $\R^d$ to general LCA groups.
We will develop  the concept of range function and the techniques of  fiberization in this general context.
The treatment will be for shift-invariant spaces following the lines of Bownik \cite{B}.
The conclusions for doubly invariant spaces will follow via the Plancherel theorem for the
Fourier transform on LCA groups.

First we will fix some notation and set our standing assumptions that will be in force for the
remainder of the manuscript.
\begin{standing}\label{standing}
We will assume throughout the next sections that.
\begin{enumerate}
\item[$\bullet$] $G$ is a second countable LCA group.
\item[$\bullet$] $H$ is a countable uniform lattice on $G$.
\end{enumerate}

We denote, as before, by $\G$  the dual group of $G$, by $\Dd$ the annihilator of $H$, and by
$\W$ a fixed Borel section of $\G/\Dd$.

The choice of particular Haar measure in each of the groups considered in this paper does not affect the validity of the results. However, different constants will appear in the formulas.

Since we have the freedom to choose the Haar measures, we will fix the following normalization in order to avoid carrying over constants through the paper and to simplify the statements of the results.

We choose $m_{\Dd}$ and $m_{H}$ such that $m_{\Dd}(\{0\})=m_{H}(\{0\})=1$.  We fix $m_{\G/\Dd}$ such that $m_{\G/\Dd}(\G/\Dd)=1$ and therefore the inversion formula holds between $H$ and  $\G/\Dd$.
Then, we set $m_{\G}$
such that Theorem \ref{teo-coset-descomposition} holds for $m_{\G}$, $m_{\G/\Dd}$ and $m_{\Dd}$.
Finally, we  normalize $m_G$ such that the inversion formula holds for $m_{\G}$ and $m_G$.
\end{standing}

As a consequence of the normalization given above and Lemma  \ref{lema-medidas},  it holds that $m_{\G}(\W)=1$.
Note that under our Standing Assumptions \ref{standing}, Theorem \ref{theo-standing} applies. So we will use
the properties of $G$, $H$, $\G$ and $\Dd$ stated in that theorem.

\subsection{Preliminaries}\noindent

The space $L^2(\W,\ell^2(\Dd))$
is the space of all measurable functions $\Phi:\W\rightarrow\ell^2(\Dd)$ such that
$$\int_{\W}\|\Phi(\w)\|_{\ell^2(\Dd)}^2\,dm_{\G}(\w)<\infty,$$
where a function $\Phi:\W\rightarrow\ell^2(\Dd)$ is measurable, if and only if for each $a\in \ell^2(\Dd)$ the function
$\w\mapsto\langle \Phi(\w),a\rangle$ is a measurable function from  $\W$
into $\C$.

\begin{remark}
This is the usual notion of weak measurability for vector functions.
If the values of the functions are in a separable space, as a consequence of Petti's Theorem,
the notions of weak and strong measurability agree.
As we have seen in Section \ref{sec-1} and according to our hypotheses,
$\Dd$ is a countable uniform lattice on $\G$.
Therefore, $\ell^2(\Dd)$ is a separable Hilbert space. Then,
in $L^2(\W,\ell^2(\Dd))$ we have only one measurability notion.
\end{remark}

The space $L^2(\W,\ell^2(\Dd))$, with the inner product
$$\langle \Phi, \Psi\rangle:=\int_{\W}\langle\Phi(\w),\Psi(\w)\rangle_{\ell^2(\Dd)}\,dm_{\G}(\w)$$
is a complex Hilbert space.

Note that for $\Phi\in L^2(\W,\ell^2(\Dd))$ and $\w\in\W$
$$\|\Phi(\w)\|_{\ell^2(\Dd)}=\Big(\sum_{\dd\in\Dd}|(\Phi(\w))_{\dd}|^2\Big)^{1/2},$$
where $(\Phi(\w))_{\dd}$ denotes the value of the sequence $\Phi(\w)$ in $\dd$.
If $\Phi\in L^2(\W,\ell^2(\Dd))$, the  sequence $\Phi(\w)$ is the {\it fiber of $\Phi$
at $\w$}.

The following proposition shows that the space $L^2(\W,\ell^2(\Dd))$ is isometric to $L^2(G)$.

\begin{proposition}\label{prop-tau}
The mapping $\T:L^2(G)\longrightarrow L^2(\W,\ell^2(\Dd))$ defined as
$$\T f(\w)=\{\widehat{f}(\w+\delta)\}_{\delta\in\Delta},$$
is an isomorphism that satisfies
$\|\T f\|_2=\|f\|_{L^2(G)}.$
\end{proposition}

The next periodization lemma will be necessary for the proof of Proposition  \ref{prop-tau}.

\begin{lemma}\label{lema-periodizar}
Let $g\in L^2(\G)$. Define the function $\mathcal{G}(\w)=\sum_{\dd\in\Dd}|g(\w+\dd)|^2$.
Then, $\mathcal{G}\in L^1(\W)$ and moreover
$$\|g\|_{L^2(\G)}=\|\mathcal{G}\|_{L^1(\W)}.$$
\end{lemma}

\begin{proof}
Since $\Omega $ is a section of  $\Gamma/\Delta$, we have
that $\Gamma=\bigcup_{\delta\in\Delta}\Omega-\delta$, where the union is disjoint. Therefore,
\begin{eqnarray*}
\int_{\Gamma}|g(\gamma)|^2\,dm_{\G}(\gamma)&=&\sum_{\delta\in\Delta}\int_{\Omega-\delta}|g(\w)|^2\,dm_{\G}(\w)\\
&=&\sum_{\delta\in\Delta}\int_{\Omega}|g(\w+\delta)|^2\,dm_{\G}(\w)\\
&=&\int_{\Omega}\sum_{\delta\in\Delta}|g(\w+\delta)|^2\,dm_{\G}(\w).\\
\end{eqnarray*}
This proves that $\mathcal{G}\in L^1(\W)$ and $\|g\|_{L^2(\G)}=\|\mathcal{G}\|_{L^1(\W)}.$
\end{proof}

\begin{proof}[Proof of Proposition \ref{prop-tau}]
First we prove that $\T$ is well defined.
For this we must show that,
$\forall f\in L^2(G)$, the vector function
$\T f$ is measurable and $\|\T f\|_2<\infty$.

According to Lemma \ref{lema-periodizar}, the sequence
$\{\widehat{f}(\w+\delta)\}_{\delta\in\Delta}\in\ell^2(\Dd)$, a.e. $\w\in\W$, for all $f\in L^2(G)$.
Then, given $a=\{a_{\dd}\}_{\dd\in\Dd}\in \ell^2(\Dd)$, the product
$\langle\T f(\w), a\rangle=\sum_{\dd\in\Dd}\widehat{f}(\w+\dd)a_{\dd}$ is finite a.e. $\w\in\W$.
From here the measurability of $f$ implies that
 $\w\mapsto \langle\T f(\w), a\rangle$ is
a measurable function in the usual sense. This proves the measurability of $\T f$.

If $f\in L^2(G)$, as a consequence of  Lemma \ref{lema-periodizar}, we have
\begin{eqnarray*}
\|\T f\|_2^2&=&\int_{\Omega}\|\T f(\w)\|^2_{\ell^2(\Delta)}\,dm_{\G}(\w) \\
&=&\int_{\Omega}\sum_{\delta\in \Delta}|\widehat{f}(\w +\delta)|^2\,dm_{\G}(\w)\\
&=&\int_{\Gamma}|\widehat{f}(\gamma)|^2\,dm_{\G}(\gamma)\\
&=&\int_G|f(x)|^2\,dm_{G}(x).\\
\end{eqnarray*}
Thus, $\|\T f\|_2<\infty$ and this also proves that $\|\T f\|_2=\|f\|_{L^2(G)}.$

What is left is to show  that  $\T$ is onto. So, given $\Phi\in L^2(\Omega,\ell^2(\Delta))$
let us see that there exists a function
$f\in L^2(G)$ such that $\T f=\Phi$.
Using that the Fourier transform is an isometric isomorphism between $L^2(G)$ and $L^2(\Gamma)$, it will be
sufficient to find
$g\in L^2(\Gamma)$ such that $\{g(\w+\delta)\} _{\delta\in\Delta}=\Phi(\w)$
a.e. $\w \in \Omega$ and then take
$f\in L^2(G)$  such that $\widehat{f}=g$.

Given $\g\in\G$, there exist unique $\w\in\W$ and $\dd\in\Dd$ such that $\g=\w+\dd$.
So, we define  $g(\g)$  as
$$g(\gamma)=\big(\Phi(\w)\big)_{\delta}.$$

The measurability of $g$ is straightforward.

Once again, according to Lemma \ref{lema-periodizar},
\begin{eqnarray*}
\int_{\Gamma}|g(\gamma)|^2\,dm_{\G}(\gamma)
&=&\int_{\Omega}\sum_{\delta\in\Delta}|g(\w+\delta)|^2\,dm_{\G}(\w)\\
&=&\int_{\Omega}\sum_{\delta\in\Delta}|\big(\Phi(\w)_{\delta}\big)|^2\,dm_{\G}(\w)\\
&=&\int_{\Omega}\|\Phi(\w)\|^2_{\ell^2(\Delta)}\,dm_{\G}(\w)\\
&=&\|\Phi\|_2^2<+\infty.\\
\end{eqnarray*}
Thus, $g\in L^2(\Gamma)$ and this completes the proof.
\end{proof}

The mapping $\T$ will be important to study the properties of functions of $L^2(G)$
in terms of their fibers, (i.e. in terms of the fibers $\T f(\w)$).

\subsection{$H$-invariant Spaces and Range Functions}\noindent

\begin{definition}
We say that a closed subspace  $V\subset L^2(G)$ is {\it $H$-invariant}
if
$$ f\in V \Rightarrow t_hf\in V\quad\forall\,\,h\in H,$$
where $t_yf(x)=f(x-y)$ denotes the translation of $f$ by an element $y$ of $G$.
\end{definition}

For a subset $\mathcal{A}\subset L^2(G)$, we define
$$E_H(\A)= \{t_h\p: \p\in\A, h\in H\}\quad\textrm{and}\quad S(\A)=\clspan \,E_H(\A).$$
We call $S(\A)$
the $H$-invariant space generated by $\A$. If $\A=\{\p\}$ ,
we simply write  $E_H(\p)$ and $S(\p)$, and we call  $S(\p)$ a {\it principal  $H$-invariant space}.

Our main  goal  is to give a characterization of $H$-invariant spaces.
We first need  to introduce the concept of range function.
\begin{definition}
A {\it range function} is a mapping $$J:\Omega\longrightarrow\{\textrm{closed subspaces of}\,\,
\ell^2(\Delta)\}.$$
\end{definition}

The subspace $J(\w)$ is called the {\it fiber space} associated to $\w.$

For a given range function $J$, we  associate to each $\w\in\W$ the orthogonal projection
onto $J(\w)$,  $P_{\w}:\ell^2(\Delta)\rightarrow J(\w)$.

A range function $J$ is {\it measurable} if for each $a\in\ell^2(\Dd)$ the function
$\w\mapsto P_{\w}a$, from $\W$ into $\ell^2(\Dd)$, is measurable.
That is,  for each $a,b\in\ell^2(\Dd)$, $\w\mapsto\langle P_{\w}a,b\rangle$ is measurable
in the usual sense.

\begin{remark}
Note that $J$ is a measurable range function if and only if for all
$\Phi\in L^2(\Omega,\ell^2(\Delta))$, the function $\w\mapsto P_{\w}\big(\Phi(\w)\big)$
is measurable. That is, $\forall\,\, b\in \ell^2(\Dd)$,
$\w\mapsto \langle P_{\w}\big(\Phi(\w)\big), b\rangle$
is measurable in the usual sense.
\label{funcion-rango-medible-equivalencia}
\end{remark}

Given a range function $J$ (not necessarily measurable) we define the  subset  $M_J$ as
$$M_J=\big\{\Phi\in L^2(\Omega,\ell^2(\Delta)): \Phi(\w)\in J(\w) \,\,a.e. \,\,\,\w\in\Omega\big\}.$$

\begin{lemma}\label{mj-cerrado}
The subset $M_J$ is closed in $L^2(\Omega,\ell^2(\Delta))$.
\end{lemma}
\begin{proof}
Let $\{\Phi_j\}_{j\in\N}\subset M_J$ such that $\Phi_j\rightarrow\Phi$ when $j\to\infty$
in $L^2(\Omega,\ell^2(\Delta))$.
Let us consider the functions  $g_j:\W\rightarrow\R_{\geq 0}$ defined as
$g_j(\w):=\|\Phi_j(\w)-\Phi(\w)\|^2_{\ell^2(\Delta)}.$
Then, $g_j$ is measurable for all $j\in\N$ and $\forall\,\,\,\alpha>0$ it holds that
$$m_{\G}(\{g_j>\f\})\leq\frac1{\f}\int_{\W}g_j(\w)\,dm_{\G}(\w)=
\frac1{\f}\int_{\W}\|\Phi_j(\w)-\Phi(\w)\|^2_{\ell^2(\Delta)}\,dm_{\G}(\w)\to 0,$$
when $j\to\infty$. So, $g_j\to 0$ in measure and therefore,
there exists a  subsequence $\{g_{j_k}\}_{k\in\N}$ of $\{g_j\}_{j\in\N}$ which goes to zero a.e. $\w\in\W$.
Then, $\Phi_{j_k}(\w)\to\Phi(\w)$ in $\ell^2(\Dd)$ a.e.\,\,$\w\in\W$ and hence, since
$\Phi_{j_k}(\w)\in J(\w)$ a.e.\,\,$\w\in\W$ and  $J(\w)$ is closed, $\Phi(\w)\in J(\w)$
a.e.\,\,$\w\in\W$. Therefore $\Phi\in M_J$.
\end{proof}

The following proposition is a generalization to the context of groups of a lemma of Helson,
(see \cite{helson} and also \cite{B}).

\begin{proposition}\label{lema-helson}
Let $J$ be a measurable range function and $P_{\w}$ the associated orthogonal projections.
Denote by $\mathcal{P}$ the orthogonal projection onto $M_J$.
Then,
$$\big(\mathcal{P}\Phi\big)(\w)=P_{\w}\big(\Phi(\w)\big),\,\,\,a.e.\,\, \w\in\W, \;\;\forall\,\, \Phi\in L^2(\Omega,\ell^2(\Delta)).$$

\end{proposition}

\begin{proof}
Let $\mathcal{Q}:L^2(\Omega,\ell^2(\Delta))\rightarrow L^2(\Omega,\ell^2(\Delta))$ be the linear mapping
$\Phi\mapsto\mathcal{Q}\Phi$, where
$$\big(\mathcal{Q}\Phi\big)(\w):=P_{\w}\big(\Phi(\w)\big).$$ We want to show that
$\mathcal{Q}=\mathcal{P}$.

Since $J$ is a measurable range function, due to Remark \ref{funcion-rango-medible-equivalencia}, $\mathcal{Q}\Phi$ is measurable
for each $\Phi\in L^2(\Omega,\ell^2(\Delta))$. Furthermore,
since $P_{\w}$ is an orthogonal projection, it has norm  one, and therefore
\begin{eqnarray*}
\|\mathcal{Q}\Phi\|_2^2&=&\int_{\W}\|\big(\mathcal{Q}\Phi\big)(\w)\|^2_{\ell^2(\Dd)}\,dm_{\G}(\w)\\
&=& \int_{\W}\|P_{\w}\big(\Phi(\w)\big)\|^2_{\ell^2(\Dd)}\,dm_{\G}(\w)\\
&\leq&\int_{\W}\|\Phi(\w)\|^2_{\ell^2(\Dd)}\,dm_{\G}(\w)=\|\Phi\|^2_2<\infty.
\end{eqnarray*}
Then,
$\mathcal{Q}$ is well defined and it has norm less or equal to 1.

From the fact that  $P_{\w}$ is an orthogonal projection,  it follows that
$\mathcal{Q}^2=\mathcal{Q}$ and $\mathcal{Q}^{*}=\mathcal{Q}$. So, $\mathcal{Q}$ is
also an orthogonal projection. To complete our proof
let us see that $M=M_J$, where $M:=Ran(\mathcal{Q})$.

By  definition of $\mathcal{Q}$, $M\subset M_J$.

If we suppose that $M$ is properly included in
$M_J$,then there exists
$\Psi\in M_J$ such that $\Psi\neq 0$ and $\Psi\perp M$.
Then, $\forall \,\,\Phi \in L^2(\Omega,\ell^2(\Delta))$,
$0=\langle\mathcal{Q}\Phi,\Psi\rangle=\langle\Phi, \mathcal{Q}\Psi\rangle$.

Hence, $\mathcal{Q}\Psi=0$ and therefore $P_{\w}\big(\Psi(\w)\big)=0$ a.e. $\w\in\W$.
Since $\Psi\in M_J$,  $\Psi(\w)\in J(\w)$  a.e. $\w\in\W$, thus
$P_{\w}\big(\Psi(\w)\big)=\Psi(\w)$  a.e. $\w\in\W$. Finally,  $\Psi=0$  a.e. $\w\in\W$
and this is a contradiction.
\end{proof}

We now give a characterization of $H$-invariant spaces using
range functions.

\begin{theorem}\label{rango-SIS}
Let $V\subset L^2(G)$ be a closed subspace and $\T$ the map defined in Proposition \ref{prop-tau}. Then,
$V$ is $H$-invariant if and only if there exists a measurable range function $J$ such that
$$V=\big\{f\in L^2(G):\,\, \T f(\w)\in J(\w)\,\,\textrm{a.e.}\,\,\w\in\W\big\}.$$

Identifying range functions which are equal almost everywhere, the correspondence between
 $H$-invariant spaces and  measurable range functions is one to one and onto.

Moreover, if $V=S(\A)$ for some countable subset $\A$ of $L^2(G)$, the measurable range function $J$
associated to $V$ is given by
$$J(\w)=\clspan\{\T\p(\w): \,\,\p\in\A\}, \,\, \text{ a.e. } \w\in\W.$$
\end{theorem}

For the proof, we need the following results.

\begin{lemma}\label{rango-1-a-1}
If $J$ and $K$ are two measurable range functions such that $M_J=M_K$, then $J(\w)=K(\w)$ a.e. $\w\in\W$.
That is, $J$ and $K$ are equal almost everywhere.
\end{lemma}

\begin{proof}
Let $P_{\w}$ and $Q_{\w}$ be the  projections associate to $J$ and $K$ respectively.
If $\mathcal{P}$ is the orthogonal projection onto $M_J=M_K$, by Proposition \ref{lema-helson} we have that, for each
$\Phi\in L^2(\Omega,\ell^2(\Delta))$
$$\big(\mathcal{P}\Phi\big)(\w)=P_{\w}\big(\Phi(\w)\big)\;\;
\textrm{ and }\;\;\big(\mathcal{P}\Phi\big)(\w)=Q_{\w}\big(\Phi(\w)\big)\,\,\,\textrm{a.e.}\,\,\w\in\W.$$

So, $P_{\w}\big(\Phi(\w)\big)=Q_{\w}\big(\Phi(\w)\big)$ a.e. $\w\in\W$, for all $\Phi\in L^2(\Omega,\ell^2(\Delta))$.
In particular, if $e_{\lambda}\in\ell^2(\Dd)$ is defined by $(e_{\lambda})_{\dd}=1$ if $\dd=\lambda$ and
$(e_{\lambda})_{\dd}=0$ otherwise, $P_{\w}(e_{\lambda})=Q_{\w}(e_{\lambda})$ a.e. $\w\in\W$, for all $\lambda\in\Dd$. Hence, since $\{e_{\lambda}\}_{\lambda\in\Dd}$ is a basis for $\ell^2(\Delta)$, it follows that
$P_{\w}=Q_{\w}$ a.e. $\w\in\W$.
Thus $J(\w)=K(\w)$ a.e. $\w\in\W$.
\end{proof}

\begin{remark}\label{tau-fourier}
Note that for  $f\in L^2(G)$ and for $h\in H$,
$$\T t_hf(\w)=(h,-\w)\T f(\w),$$
since $\forall\,\,y\in G$, $\widehat{t_yf}(\g)=(y,-\g)\widehat{f}(\g)$ and, as we showed in Remark \ref{caracteres-delta-periodicos}, the character $(h,.)$ is
$\Dd$-periodic.
\end{remark}

\begin{proof}[Proof of theorem \ref{rango-SIS}]
Let us first suppose that $V$ is $H$-invariant.
Since $L^2(G)$ is separable,
$V=S(\A)$ for some countable subset $\A$ of  $L^2(G)$.

We define the  function $J$ as $J(\w)=\clspan\{ \T\p(\w): \,\,\p\in\A\}$.
Note that since $\A$ is a countable set,
 $J$ is well defined a.e. $\w\in\W$. We will prove that $J$ satisfies:
\begin{enumerate}
\item[(i)] $V=\big\{f\in L^2(G):\,\, \T f(\w)\in J(\w)\,\,\textrm{a.e.}\,\,\w\in\W\big\}$,
\item[(ii)] $J$ is measurable.
\end{enumerate}

To show (i) it is sufficient to prove that  $M=M_J$, where $M:=\T V$.
Let $\Phi\in M$. Then, $\T^{-1}\Phi\in V=\clspan\{t_h\p:\,\,h\in H,\,\,\p\in\A\}$. Therefore,
there exists a sequence $\{g_j\}_{j\in\N}\subset \Span \{t_h\p:\,\,h\in H,\,\,\p\in\A\}$
such that $\T g_j:=\Phi_j$ converges in  $L^2(\W,\ell^2(\Dd))$ to $\Phi$, when $j\to\infty$ .

Due to the  definition of $J$ and  Remark \ref{tau-fourier},
$\Phi_j(\w)\in J(\w)$ a.e. $\w\in\W$.
Thus, in the same way that in Lemma \ref{mj-cerrado}, we can prove that $\Phi(\w)\in J(\w)$
a.e. $\w\in\W$ and therefore $\Phi\in M_J$. So, $M\subset M_J$.

Let us suppose that there exists $\Psi\in L^2(\W,\ell^2(\Dd))$, such that $\Psi \neq 0$ and $\Psi$ is orthogonal to $M$.
Then, for each $\Phi\in M$, $\langle\Phi,\Psi\rangle=0.$
In particular, if $\Phi\in\T\A\subset\T V=M$ and $h\in H$, we have that
$(h,.)\Phi(.)\in \T V=M$ since $(h,.)\Phi(.)=\T(t_{-h}\T^{-1}\Phi)(.)$ and $t_{-h}\T^{-1}\Phi\in V$.

So, as  $(h,.)$ is $\Dd$-periodic,
$$0=\langle (h,.)\Phi(.),\Psi\rangle
=\int_{\W}(h,\w)\langle \Phi(\w),\Psi(\w)\rangle_{\ell^2(\Dd)}dm_{\G}(\w).$$

Hence, by Proposition \ref{prop-fourier-H-0}, $\langle \Phi(\w),\Psi(\w)\rangle_{\ell^2(\Dd)}=0$ a.e. $\w\in\W$, and this holds $\forall\,\,\Phi\in\T(\A)$.
Therefore  $\Psi(\w)\in J(\w)^{\perp}$ a.e. $\w\in\W$.

Now, if  $M$ is properly included  in $M_J$, there  exists $\Psi\in M_J$, with $\Psi \neq 0$ and orthogonal
to $M$.
Hence, $\Psi(\w)\in J(\w)^{\perp}$ a.e. $\w\in\W$. On the other hand since $\Psi\in M_J$,
$\Psi(\w)\in J(\w)$ a.e. $\w\in\W$.
Thus,  $\Psi(\w)=0$ a.e. $\w\in\W$ and this is a contradiction.
Therefore $M=M_J$.

It remains to prove  that the  range function $J$ is measurable.
For this  we must show that, for all $a,b\in\ell^2(\Dd)$, $\w\mapsto\langle P_{\w}a,b\rangle$ is measurable, where $P_{\w}:\ell^2(\Dd)\rightarrow J(\w)$ are the orthogonal projections associated to $J(\w)$,

Let $\mathcal{I}$ be the identity mapping in  $L^2(\W,\ell^2(\Dd))$ and $\mathcal{P}:L^2(\W,\ell^2(\Dd))\rightarrow M$ the orthogonal projection associated to $M$.
If $\Psi\in L^2(\W,\ell^2(\Dd))$, the function $(\mathcal{I}-\mathcal{P})\Psi$ is  orthogonal to $M$ and,
by the above reasoning,
$(\mathcal{I}-\mathcal{P})\Psi(\w)\in J(\w)^{\perp}$, a.e. $\w\in\W$.
Then,
$$P_{\w}((\mathcal{I}-\mathcal{P})\Psi(\w))=P_{\w}(\Psi(\w)-\mathcal{P}\Psi(\w))=0$$
 a.e. $\w\in\W$ and therefore  $P_{\w}(\Psi(\w))=P_{\w}(\mathcal{P}\Psi(\w))=\mathcal{P}\Psi(\w)$
a.e $\w\in\W$.
In particular, $P_{\w}a=\mathcal{P}a(\w)$ a.e. $\w\in\W$, $\forall\,\,\, a\in \ell^2(\Dd)$. Thus, since
$\w\mapsto\langle \mathcal{P}a(\w)a,b\rangle$ is measurable $\forall\,\,b\in\ell^2(\Dd)$,
$\w\mapsto\langle P_{\w}a,b\rangle$ is measurable as well.

Conversely, if $J$ is a measurable range function, let us see that the closed  subspace in $L^2(G)$,
defined by $V:=\T^{-1}(M_J)$ is $H$-invariant.
For this, let us consider  $f\in V$ and $h\in H$ and  let us prove that $t_hf\in V$.

Since
$\T(t_hf)(\w)=(h,-\w)\T f(\w)$ a.e. $\w\in\W$ and  $\T f\in M_J$,
we have that  $(h,-\w)\T f(\w)\in J(\w)$ a.e. $\w\in\W$.
Then, $\T(t_hf)\in M_J$ and  therefore $t_hf\in V$.

Furthermore, $V=S(\A)$ for some  countable set $\A$ of $L^2(G)$. Then,
$$K(\w)=\clspan\{\T\p(\w): \,\,\p\in\A\}\,\,\,\textrm{a.e.}\,\,\w\in\W,$$
defines a measurable range function which satisfies  $V=\T^{-1}(M_K)$.
Thus, $M_K=\T V=M_J$. Since $J$ and $K$ are both measurable range functions, Lemma \ref{rango-1-a-1} implies that $J=K$ a.e $\w\in\W$.

This also shows that the correspondence between $V$ and $J$ is onto and one to one.
\end{proof}

\section{Frames and Riesz Basis for $H$-invariant Spaces}\label{sec-3}
Let $\mathcal{H}$ be a Hilbert space and $\{u_i\}_{i\in I}$ a sequence of  $\mathcal{H}$.

The sequence $\{u_i\}_{i\in I}$ is a {\it Bessel sequence} in $\mathcal{H}$ with constant $B$ if
$$\sum_{i\in I}|\langle f, u_i\rangle|^2\leq B\|f\|^2,\;\;\textrm{for all}\;\; f\in\mathcal{H}.$$
The sequence $\{u_i\}_{i\in I}$ is a {\it frame} for $\mathcal{H}$ with constants $A$ and $B$ if
$$A\|f\|^2\leq \sum_{i\in I}|\langle f, u_i\rangle|^2\leq B\|f\|^2,\;\;\textrm{for all}\;\; f\in\mathcal{H}.$$
The frame $\{u_i\}_{i\in I}$ is a {\it tight frame} if $A=B$,
and the frame $\{u_i\}_{i\in I}$ is a {\it Parseval frame} if $A=B=1$.

The sequence $\{u_i\}_{i\in I}$ is a {\it Riesz sequence} for $\mathcal{H}$ if there exist
positive constants $A$ and $B$  such that
$$A\sum_{i\in I}|a_i|^2\leq||\sum_{i\in I}a_iu_i||^2_{\mathcal{H}}\leq B \sum_{i\in I}|a_i|^2$$
for all $\{a_i\}_{i\in I}$ with finite support.
Moreover, a Riesz sequence that in addition, is  a complete family in  $\mathcal{H}$, is a  {\it Riesz basis} for $\mathcal{H}$.

We are now ready to prove a result which  characterizes when $E_H(\A)$ is a frame of
$L^2(G)$ in terms of the fibers $\{\T\p(\w)\,:\,\p\in\A\}$. It generalizes
Theorem 2.3 of \cite{B} to the  context of groups.

\begin{theorem}
Let $\A$ be a countable subset of $L^2(G)$, $J$ the measurable range function associated to
$S(\A)$ and  $A\leq B$ positive constants.
Then, the following propositions are equivalent:
\begin{enumerate}
\item[(i)] \label{item-1-marcos}The set $E_H(\A)$ is a frame for $S(\A)$ with  constants
$A$ and $B$.
\item[(ii)] \label{item-2-marcos}For almost every $\w\in\W$, the set $\{\T\p(\w):\p\in\A\}\subset\ell^2(\Dd)$
is a frame for $J(\w)$ with  constants $A$ and $B$.
\end{enumerate}
\label{teo-fibras-marcos}
\end{theorem}

\begin{proof}
Since $\langle f,g\rangle_{L^2(G)}=\langle \T f,\T g\rangle_{L^2(\W,\ell^2(\Dd))}$,
by Remark \ref{tau-fourier} we have that
\begin{equation*}
\begin{split}
\sum_{h\in H}\sum_{\p\in\A}|\langle t_h\p,f&\rangle_{L^2(G)}|^2
=\sum_{h\in H}\sum_{\p\in\A}|\langle \T(t_h\p),\T f\rangle_{L^2(\W,\ell^2(\Dd))}|^2\\
&=\sum_{\p\in\A}\sum_{h\in H}
|\int_{\W}(h,-\w)\langle \T\p(\w),\T f(\w)\rangle_{\ell^2(\Dd)}dm_{\G}(\w)|^2.\\
\end{split}
\end{equation*}

Let us define for each $ \varphi \in \A$, the following,
$$ R(\p) = \sum_{h\in H}|\int_{\W}(h,-\w)\langle \T\p(\w),\T f(\w)\rangle_{\ell^2(\Dd)}dm_{\G}(\w)|^2$$
and
$$ T(\p) = \int_{\W}|\langle \T\p(\w),\T f(\w)\rangle_{\ell^2(\Dd)}|^2dm_{\G}(\w).$$

$(i)$ $\Rightarrow$ $(ii)$
If $E_H(\A)$ is a frame for $S(\A)$, in particular it holds that $\forall f\in S(\A)$,
$\sum_{h\in H\,\p\in\A}|\langle t_h\p, f\rangle|^2<\infty.$

Then, for each $\p\in\A$, we have that $R(\p)<\infty.$
Therefore, the sequence $\{c_h\}_{h\in H}$, with
$$c_h:=\int_{\W}(h,\w)\langle \T\p(\w),\T f(\w)\rangle_{\ell^2(\Dd)}dm_{\G}(\w),$$
belongs to $\ell^2(H)$.

Let us consider the function $F(\w):=\sum_{h\in H}c_h\eta_h(\w)$, where
$\eta_h$ are the functions defined in Lemma \ref{lema-bon}.
Then, since $\{c_h\}_{h\in H}\in \ell^2(H)$ and $\{\eta_h\}_{h\in H}$
is an orthogonal basis of $L^2(\W)$,
we have that $F\in L^2(\W)\subset L^1(\W)$ (recall that $m_{\G}(\W)<\infty$).

On the other hand,
the function $\psi(\w):=\langle \T\p(\w),\T f(\w)\rangle_{\ell^2(\Dd)}$ belongs to $L^1(\W)$.
So,
$\psi-F\in L^1(\W)$ and moreover
$$\int_{\W}(h,-\w)(\psi(\w)-F(\w))dm_{\G}(\w)= c_{-h}-c_{-h}=0$$ for all $h\in H$.
Thus, Proposition \ref{prop-fourier-H-0}
yields that $F=\psi$ a.e. $\w\in\W$.
Therefore $\psi\in L^2(\W)$ and
$$\psi(\w)=\sum_{h\in H}c_h\eta_h(\w),$$
a.e. $\w\in\W$.

As a consequence of Proposition \ref{lema-parseval}, we obtain that $R(\p) = T(\p)$ holds for all $\p\in\A$.

We will now prove  that, for almost every  $\w\in\W$, $\{\T\p(\w):\p\in\A\}$
is a frame with  constants $A$ and $B$
for $J(\w)$.

Let us suppose that
\begin{equation}
A\|P_{\w}d\|^2_{\ell^2(\Dd)}\leq\sum_{\p\in\A}|\langle\T\p(\w), P_{\w}d\rangle|^2
\leq B\|P_{\w}d\|^2_{\ell^2(\Dd)}
\label{ecu-1-marcos}
\end{equation}
a.e. $\w\in\W$, for each $d\in\mathcal{D}$,
where  $\mathcal{D}$ is a dense countable subset of $\ell^2(\Dd)$ and $P_{\w}$
are the orthogonal projections associated to $J$.
Then, for each $d\in\mathcal{D}$, let $Z_d\subset\W$ be a measurable set with $m_{\G}(Z_d)=0$ such that (\ref{ecu-1-marcos})
holds for all $\w\in \W\setminus Z_d$.
So the set $Z=\bigcup_{d\in\mathcal{D}}Z_d$ has null $m_{\G}$-measure. Therefore for $\w\in \W\setminus Z$ and $a\in J(\w)$,
using  a density argument it follows from (\ref{ecu-1-marcos}) that
$$A\|a\|^2\leq\sum_{\p\in\A}|\langle\T\p(\w), a\rangle|^2\leq B\|a\|^2.$$

Thus, it is sufficient  to show that (\ref{ecu-1-marcos}) holds.
For this, we will suppose that this is not so and we will prove that there exist
$d_0\in\mathcal{D}$, a measurable set $W\subset\W$ with $m_{\G}(W)>0$, and $\varepsilon>0$ such that
$$\sum_{\p\in\A}|\langle\T\p(\w), P_{\w}d_0\rangle|^2>(B+\e)\|P_{\w}d_0\|^2,\,\,\forall\,\,\w\in W$$
or
$$\sum_{\p\in\A}|\langle\T\p(\w), P_{\w}d_0\rangle|^2<(A-\e)\|P_{\w}d_0\|^2,\,\,\forall\,\,\w\in W.$$

So,  let us take $d_0\in\mathcal{D}$ for which  (\ref{ecu-1-marcos}) fails.
Then at least one of this sets
$$\{\w\in\W : K(\w)-B\|P_{\w}d_0\|^2>0\}\quad,\quad
\{\w\in\W : K(\w)-A\|P_{\w}d_0\|^2<0\}$$
has positive measure, where $K(\w):=\sum_{\p\in\A}|\langle\T\p(\w), P_{\w}d_0\rangle|^2$.
Let us suppose, without loss of generality, that $$m_{\G}(\{\w\in\W : K(\w)-B\|P_{\w}d_0\|^2>0\})>0.$$
Since
\begin{equation*}
\{\w\in\W : K(\w)-B\|P_{\w}d_0\|^2>0\}
=\bigcup_{j\in\N}\{\w\in\W : K(\w)-B+\frac1{j})\|P_{\w}d_0\|^2>0\},
\end{equation*}
there exists at least one set in the union, in the right hand side  of this equality, with positive
measure and this proves our claim.

Then, we can suppose that
\begin{equation}
\sum_{\p\in\A}|\langle\T\p(\w), P_{\w}d_0\rangle|^2>(B+\e)\|P_{\w}d_0\|^2,\,\,\forall\,\,\w\in W
\label{ecu-marco-falla}
\end{equation}
holds.
Now take $f\in S(\A)$ such that $\T f(\w)=\chi_{W}(\w)P_{\w}d_0$.
Note that this is  possible since, by Theorem \ref{rango-SIS}, $\chi_{E}(\w)P_{\w}d_0$ is a measurable function.

As $E_H(\A)$ is a frame for $S(\A)$ and
$$\sum_{h\in H}\sum_{\p\in\A}|\langle t_h\p, f\rangle_{L^2(G)}|^2
=\sum_{\p\in\A}\int_{\W}|
\langle \T\p(\w),\T f(\w)\rangle_{\ell^2(\Dd)}|^2dm_{\G}(\w),$$
we have that
\begin{equation}\label{ecu-5-marcos}
A\|f\|^2\leq
\sum_{\p\in\A}\int_{\W}|\langle\T\p(\w), \T f(\w)\rangle|^2dm_{\G}(\w)\leq B\|f\|^2.
\end{equation}

Using Proposition \ref{prop-tau}, we can rewrite (\ref{ecu-5-marcos}) as
\begin{equation}
A\|\T f\|^2\leq
\sum_{\p\in\A}\int_{\W}|\langle\T\p(\w), \T f(\w)\rangle|^2dm_{\G}(\w)\leq B\|\T f\|^2.
\label{ecu-2-marcos}
\end{equation}

Now, $$\|\T f\|^2=\int_{\W}\chi_{W}(\w)\|P_{\w}d_0\|^2dm_{\G}(\w)$$
and if we integrate in (\ref{ecu-marco-falla}) over $W$, we obtain
\begin{equation*}
\sum_{\p\in\A}\int_{\W}|\langle \T\p(\w),\chi_{W}(\w)P_{\w}d_0\rangle|^2\,dm_{\G}(\w)
\geq(B+\e)\|\T f\|^2.
\end{equation*}
This is a contradiction with  inequality (\ref{ecu-2-marcos}).
Therefore, we proved inequality (\ref{ecu-1-marcos}).

$(ii)$ $\Rightarrow$ $(i)$
If now $\{\T\p(\w):\p\in\A\}$ is a frame for $J(\w)$ a.e $\w\in\W$ with
constants $A$ and $B$, we have that
\begin{equation*}
A\|a\|^2\leq\sum_{\p\in\A}|\langle\T\p(\w), a\rangle|^2\leq B\|a\|^2
\end{equation*}
for all $a\in J(\w)$.
In particular, if $f\in S(\A)$, by Theorem \ref{rango-SIS}, $\T f(\w)\in J(\w)$ a.e. $\w\in\W$ and then,
\begin{equation}
A\|\T f(\w)\|^2\leq\sum_{\p\in\A}|\langle\T\p(\w), \T f(\w)\rangle|^2\leq B\|\T f(\w)\|^2
\label{ecu-3-marcos}
\end{equation}
a.e. $\w\in\W$.

Thus, integrating  (\ref{ecu-3-marcos}) over $\W$, we obtain
\begin{equation}
A\|\T f\|^2\leq\int_{\W}\sum_{\p\in\A}|\langle\T\p(\w), \T f(\w)
\rangle|^2dm_{\G}(\w)\leq B\|\T f\|^2.
\label{ecu-4-marcos}
\end{equation}
So, $\langle \T \p(.),\T f(.)\rangle$  belongs to $L^2(\W)$  for each $\p\in\A$
and the
equality $R(\p)= T(\p)$, can be obtained in a similar way as we did before.

Finally, since
$\|\T f\|_2^2=\|f\|_2^2$
and
$$\sum_{h\in H}\sum_{\p\in\A}|\langle t_h\p, f\rangle_{L^2(G)}|^2
=\sum_{\p\in\A}\int_{\W}|
\langle \T\p(\w),\T f(\w)\rangle_{\ell^2(\Dd)}|^2dm_{\G}(\w),$$
 inequality (\ref{ecu-4-marcos}) implies that  $E_H(\A)$ is a frame for $S(\A)$ with constants $A$ and $B$.
\end{proof}

Theorem \ref{teo-fibras-marcos} reduces the problem of  when $E_H(\A)$ is a frame for $S(\A)$
to  when the fibers $\{\T\p(\w):\p\in\A\}$  form  a frame for $J(\w)$.
The advantage of this reduction is that, for example, when $\A$ is a finite set,
the fiber spaces $\{\T\p(\w):\p\in\A\}$ are finite dimensional while  $S(\A)$ has infinite dimension.

If $\A=\{\p\}$, Theorem \ref{teo-fibras-marcos}
generalizes  a known result for the case $G=\R^d$ to the context of groups.
This is stated in the next corollary, which was  proved in \cite{KR}.
We give here a different proof.

\begin{corollary}\label{coro-ppal-marcos}
Let $\p\in L^2(\W)$ and  $\W_{\p}=\{\w\in\W\,:\,\sum_{\dd\in\Dd}|\widehat{\p}(\w+\dd)|^2\neq0\}$. Then, the following are equivalent:
\begin{enumerate}
\item[(i)] The set $E_H(\p)$ is a frame for $S(\p)$ with constants $A$ and $B$.
\item[(ii)] $A\leq \sum_{\dd\in\Dd}|\widehat{\p}(\w+\dd)|^2
\leq B,$ a.e. $\w\in\W_{\p}$.
\end{enumerate}
\end{corollary}

\begin{proof}
Let $J$ be the measurable range function associated to $S(\p)$.
Then, by Theorem \ref{rango-SIS}, $J(\w)=\textrm{span}\{\T\p(\w)\}$ a.e $\w\in\W$.
Thus, each $a\in J(\w)$ can be written as $a=\lambda\T\p(\w)$ for some $\lambda\in\C$.

Therefore, by
Theorem \ref{teo-fibras-marcos},  $(i)$ holds if and only if,  for almost every $\w\in\W$ and for all $\lambda\in\C$,
\begin{equation}\label{ec-coro-marco}
A\|\lambda\T\p(\w)\|^2\leq|\lambda|^2\|\T\p(\w)\|^4\leq B\|\lambda\T\p(\w)\|^2.
\end{equation}

Then, since $\|\T\p(\w)\|^2=\sum_{\dd\in\Dd}|\widehat{\p}(\w+\dd)|^2$, (\ref{ec-coro-marco})
holds if and only if
$$A\leq \sum_{\dd\in\Dd}|\widehat{\p}(\w+\dd)|^2
\leq B,\quad {\text  a.e. } \,\w\in\W_{\p}.$$
\end{proof}

For the case of Riesz basis, we have an analogue result to Theorem \ref{teo-fibras-marcos}.

\begin{theorem}\label{teo-fibras-bases-riesz}
Let $\A$ be a countable subset of $L^2(G)$, $J$ the measurable range function associated to
$S(\A)$ and  $A\leq B$ positive constants.
Then, they are equivalent:
\begin{enumerate}
\item[(i)] The set $E_H(\A)$ is a  Riesz basis for $S(\A)$ with
constants $A$ and $B$.
\item[(ii)] For almost every $\w\in\W$, the set $\{\T\p(\w):\p\in\A\}\subset\ell^2(\Dd)$
is a Riesz basis for $J(\w)$ with constants $A$ and $B$.
\end{enumerate}
\end{theorem}

For the proof we will need the next lemma.

\begin{lemma}
For each $m\in L^{\infty}(\W)$ there exists a sequence of trigonometric polynomials $\{P_k\}_{k\in\N}$
such that:
\begin{enumerate}
\item[(i)] $P_k(\w)\rightarrow m(\w)$, a.e. $\w\in\W$,
\item[(ii)] There exists $C>0$, such that $\|P_k\|_{\infty}\leq C$, for all $k\in\N$.
\end{enumerate}
\label{lema-lusin}
\end{lemma}

\begin{proof}
By Lemma \ref{poly-triginometricos-densos}, taking into account Remark \ref{estructura-seccion}, we have that the trigonometric polynomials are dense in $\mathcal{C}(\W)$.

By Lusin's Theorem, for each $k\in\N$, there exists a closed set
$E_k\subset\W$ such that $m_{\G}(\W\setminus E_k)<2^{-k}$ and $m\vert_{E_k}$ is a continuous function
where $m\vert_{E_k}$ denotes the
function $m$ restricted to $E_k$.

Since $\W$ is compact, $E_k$ is compact as well.
Therefore, $m\vert_{E_k}$ is bounded.

Let $m_1, m_2:E_k\rightarrow\R$ be continuous function such that $m\vert_{E_k}=m_1+im_2$.
As a consequence of Tietze's Extension Theorem, it is possible to extend $m_1$ and $m_2$, continuously
  to all $\W$ keeping their
norms in $L^{\infty}(E_k)$. Let us call the extensions $\overline{m_1}$ and $\overline{m_2}$
and let $\overline{m_k}=\overline{m_1}+i\overline{m_2}$. Then, we have:
\begin{enumerate}
\item[(1)] $\overline{m_k}\vert_{E_k}=m\vert_{E_k}$,
\item[(2)] $\|\overline{m_k}\|_{\infty}\leq \|\overline{m_1}\|_{\infty}
+\|\overline{m_2}\|_{\infty}\leq \|m_1\|_{\infty}+\|m_2\|_{\infty}\leq 2\|m\|_{\infty}$.
\end{enumerate}

Now, by Lemma \ref{poly-triginometricos-densos}, there exists a trigonometric polynomial $P_k$
such that
$\|P_k-\overline{m_k}\|_{\infty}<2^{-k}$. So,
\begin{enumerate}
\item[({\it a})] $|P_k(\w)-m(\w)|<2^{-k}$, for all $\w\in E_k$,
\item[({\it b})] $\|P_k\|_{\infty}\leq \|P_k-\overline{m_k}\|_{\infty}
+\|\overline{m_k}\|_{\infty}\leq 2^{-k}+2\|m\|_{\infty}\leq 1+2\|m\|_{\infty}$.
\end{enumerate}

Repeating  this argument for each $k\in\N$, we obtain a  sequence $\{P_k\}_{k\in\N}$ of trigonometric polynomials   and a
 sequence  $\{E_k\}_{k\in\N}$ of sets, which satisfy  conditions $({\it a})$ and $({\it b})$.

Let $E=\cup_{j=1}^{\infty}\cap_{k=j}^{\infty}E_k.$ It is a straightforward to see that $m_{\G}(\W\setminus E)=0$.
Let us prove that if $\w\in E$, $P_k(\w)\rightarrow m(\w)$, for $k\to\infty$.
Since $\w\in E$, there exists $k_0\in\N$ for which $\w\in E_k,\,\,\,\forall \,\,k\geq k_0$.
Then, for all $k\geq k_0$, we obtain that
$|P_k(\w)-m(\w)|=|P_k(\w)-m_k(\w)|<2^{-k}\rightarrow 0$, when $k\to\infty$.
This proves part $(i)$ of this lemma and taking $C:=1+2\|m\|_{\infty}$ we have that $(ii)$ holds.
\end{proof}

\begin{proof}[Proof of Theorem \ref{teo-fibras-bases-riesz}]
Since $S(\A)=\clspan\, E_H(\A)$ and, by Theorem \ref{rango-SIS}, $J(\w)=\clspan\{\T\p(\w)\,:\,\p\in\A\}$,
we only need to show that $E_H(\A)$ is a Riesz sequence for $S(\A)$ with constants $A$ and $B$ if and only if
for almost every $\w\in\W$, the set $\{\T\p(\w):\p\in\A\}\subset\ell^2(\Dd)$
is a Riesz sequence for $J(\w)$ with constants $A$ and $B$.

For the proof of the equivalence in the theorem, we will use the
following reasoning.

Let $\{a_{\p,h}\}_{(\p,h)\in\A\times H}$ be a  sequence of finite support and let $P_{\p}$ be the
trigonometric polynomials defined  by
$$P_{\p}(\w)=\sum_{h\in H}a_{\p,h}\eta_h(\w),$$
with $\w\in\W$ and $\eta_h$ as in Proposition \ref{lema-bon}.

Note that, since $\{a_{\p,h}\}_{(\p,h)\in\A\times H}$ has finite support, only a finite number of
the  polynomials $P_{\p}$ are not  zero.

Now, as a consequence of Proposition \ref{prop-tau} we have
\begin{equation}
\begin{split}
\|\sum_{(\p,h)\in\A\times H}a_{\p,h}t_h\p\|&^2_{L^2(G)}=
\|\sum_{(\p,h)\in\A\times H}a_{\p,h}\T t_h\p\|^2_{L^2(\W,\ell^2(\Dd))}\\
&=\int_{\W}\|\sum_{(\p,h)\in\A\times H}a_{\p,h}
(-h,\w)\T \p(\w)\|^2_{\ell^2(\Dd)}\,dm_{\G}(\w)\\
&=\int_{\W}\|\sum_{\p\in\A}P_{\p}(\w)\T \p(\w)\|^2_{\ell^2(\Dd)}\,dm_{\G}(\w).\\
\end{split}
\label{ecu-riesz-1}
\end{equation}
Furthermore, by  Lemma \ref{lema-parseval},
$$\sum_{h\in H}|a_{\p,h}|^2=\|\{a_{\p,h}\}_{h\in H}\|^2_{\ell^2(H)}=
\|P_{\p}\|^2_{L^2(\W)},$$
and  adding  over $\A$, we obtain
\begin{equation}
\sum_{(\p,h)\in\A\times H}|a_{\p,h}|^2=
\sum_{\p\in\A}\|P_{\p}\|^2_{L^2(\W)}.
\label{ecu-riesz-2}
\end{equation}

$(ii)$ $\Rightarrow$ $(i)$
If we suppose that for almost every $\w\in\W$, $\{\T\p(\w):\p\in\A\}\subset\ell^2(\Dd)$
is a Riesz sequence for  $J(\w)$  with  constants $A$ and $B$,
\begin{equation}\label{ecu-riesz-0}
A\sum_{\p\in\A}|a_{\p}|^2\leq \|\sum_{\p\in\A}a_{\p}\T \p(\w)\|^2_{\ell^2(\Dd)}
\leq B\sum_{\p\in\A}|a_{\p}|^2
\end{equation}
for all $\{a_{\p}\}_{\p\in\A}$ with finite support.

In particular, the above inequality holds for $\{a_{\p}\}_{\p\in\A}=\{P_{\p}(\w)\}_{\p\in\A}$.
Now, in (\ref{ecu-riesz-0}), we can integrate over $\W$  with $\{a_{\p}\}_{\p\in\A}=\{P_{\p}(\w)\}_{\p\in\A}$, in order to obtain
\begin{equation}\label{ecu-riesz-7}
\begin{split}
A\sum_{\p\in\A}\|P_{\p}\|^2_{L^2(\W)}&\leq \int_{\W}\|\sum_{\p\in\A}P_{\p}(\w)\T
 \p(\w)\|^2_{\ell^2(\Dd)}\,dm_{\G}(\w)\\
&\leq B\sum_{\p\in\A}\|P_{\p}\|^2_{L^2(\W)}.
\end{split}
\end{equation}
Using equations (\ref{ecu-riesz-1}) and (\ref{ecu-riesz-2}) we can rewrite (\ref{ecu-riesz-7}) as
\begin{equation*}
A\sum_{(\p,h)\in\A\times H}|a_{\p,h}|^2\leq
\|\sum_{(\p,h)\in\A\times H}a_{\p,h}t_h\p\|^2_{L^2(G)}
\leq B\sum_{(\p,h)\in\A\times H}|a_{\p,h}|^2.
\end{equation*}
Therefore $E_H(\A)$ is a Riesz sequence of  $S(\A)$ with
constants $A$ and  $B$.

$(i)$ $\Rightarrow$ $(ii)$
We want to prove that, for every $a=\{a_{\p}\}_{\p\in\A}\in\ell^2(\A)$ with finite support, we have  a.e. $\w\in\W$
\begin{equation}
A\sum_{\p\in\A}|a_{\p}|^2\leq
\|\sum_{\p\in\A}a_{\p}\T \p(\w)\|^2_{\ell^2(\Dd)}
\leq B\sum_{\p\in\A}|a_{\p}|^2.
\label{ecu-riesz-3}
\end{equation}

Let us suppose that (\ref{ecu-riesz-3}) fails. Then, using a similar argument as in Theorem
\ref{teo-fibras-marcos}, we can see that there exist
$a=\{a_{\p}\}_{\p\in\A}\in\ell^2(\A)$ with finite support, a measurable set $W\subset\W$  with $m_{\G}(W)>0$
and $\varepsilon>0$ such that
\begin{equation}
\|\sum_{\p\in\A}a_{\p}\T \p(\w)\|^2_{\ell^2(\Dd)}>(B+\e)\sum_{\p\in\A}|a_{\p}|^2
,\,\,\forall\,\,\w\in W
\label{ecu-riesz-4}
\end{equation}
or
\begin{equation}
\|\sum_{\p\in\A}a_{\p}\T \p(\w)\|^2_{\ell^2(\Dd)}<(A-\e)\sum_{\p\in\A}|a_{\p}|^2
,\,\,\forall\,\,\w\in W.
\label{ecu-riesz-5}
\end{equation}

With $a=\{a_{\p}\}_{\p\in\A}$ and $W$, we define for each $\p\in\A$, $m_{\p}:=a_{\p}\chi_W$.
Thus, $m_{\p}\in L^{\infty}(\W)$ and only finitely many of these functions are not null.

By  Lemma \ref{lema-lusin}, for each $\p\in\A$ there exists a  polynomial sequence
$\{P^{\p}_k\}_{k\in\N}$ such that
\begin{enumerate}
\item[(i)] $P^{\p}_k\rightarrow m_{\p}$,
\item[(ii)] $\|P^{\p}_k\|_{\infty}\leq 1+2\|m_{\p}\|_{\infty}, \,\,\forall\,\,k\in\N$.
\end{enumerate}

Since $E_H(\A)$ is a Riesz sequence for  $S(\A)$ with
constants $A$ and  $B$,
\begin{eqnarray*}
A\sum_{(\p,h)\in\A\times H}|a_{\p,h}|^2&\leq&
\|\sum_{(\p,h)\in\A\times H}a_{\p,h}t_h\p\|^2_{L^2(G)}\\
&\leq& B\sum_{(\p,h)\in\A\times H}|a_{\p,h}|^2,
\end{eqnarray*}
for each sequence $\{a_{\p,h}\}_{(\p,h)\in\A\times H}$ with finite support.

Now, for each $k\in\N$ take  $\{a_{\p,h}\}_{(\p,h)\in\A\times H}$ to be the sequence formed with the coefficients of the polynomials  $\{P^{\p}_k\}_{\p\in\A}$.

Then, using  (\ref{ecu-riesz-1}) and (\ref{ecu-riesz-2}), we have for each $k\in\N$
\begin{equation}
A\sum_{\p\in\A}\|P^{\p}_k\|^2_{L^2(\W)}\leq \int_{\W}\|\sum_{\p\in\A}P^{\p}_k(\w)\T
\p(\w)\|^2_{\ell^2(\Dd)}\,dm_{\G}(\w)
\leq B\sum_{\p\in\A}\|P^{\p}_k\|^2_{L^2(\W)}.
\label{ecu-riesz-8}
\end{equation}
Therefore, since  $m_{\G}(\W)<\infty$ and by the Dominated Convergence Theorem,
inequality (\ref{ecu-riesz-8}) can be extended to $m_{\p}$ as
\begin{equation}
A\sum_{\p\in\A}\|m_{\p}\|^2_{L^2(\W)}\leq \int_{\W}\|\sum_{\p\in\A}m_{\p}(\w)\T
\p(\w)\|^2_{\ell^2(\Dd)}\,dm_{\G}(\w)
\leq B\sum_{\p\in\A}\|m_{\p}\|^2_{L^2(\W)}.
\label{ecu-riesz-6}
\end{equation}
So, if  (\ref{ecu-riesz-4}) occurs,  integrating over $\W$ we obtain
$$\int_{\W}\|\sum_{\p\in\A}m_{\p}(\w)\T \p(\w)\|^2_{\ell^2(\Dd)}|^2\,dm_{\G}(\w)
>(B+\e)\int_{\W}\sum_{\p\in\A}|m_{\p}(\w)|^2\,dm_{\G},$$
which  contradicts  inequality (\ref{ecu-riesz-6}). We can  proceed analogously if
(\ref{ecu-riesz-5}) occurs. Hence, (\ref{ecu-riesz-3}) holds.

\end{proof}

For the case of principal $H$-invariant spaces we have the following corollary.

\begin{corollary}\label{coro-ppal-riesz}
Let $\p\in L^2(\W)$. Then, the following are equivalent:
\begin{enumerate}
\item[(i)]  The set $E_H(\p)$ is a Riesz basis for $S(\p)$ with constants $A$ and $B$.
\item[(ii)] $A\leq \sum_{\dd\in\Dd}|\widehat{\p}(\w+\dd)|^2
\leq B,$ a.e. $\w\in\W$.
\end{enumerate}
\end{corollary}

\begin{proof}
The proof is a straightforward consequence of Theorem \ref{teo-fibras-bases-riesz} and Theorem \ref{rango-SIS}.
\end{proof}

We now  want to give another characterization of when the set $E_H(\A)$ is a frame (Riesz sequence) for $S(\A)$ with constants $A$ and $B$. For this
we will introduce what in the classical case are the synthesis and analysis operators.

For an LCA group $G$ and for a subgroup $H$ as in (\ref{standing})
let us consider a  subset
$\A=\{\p_i:i\in I\}$ of $L^2(G)$
where $I$ is a countable set.

Let $\W$ be a Borel  section of  $\G/\Dd$.
Fix $\w\in\W$ and let $\mathcal{D}$ be the set of sequences in $\ell^2(I)$ with finite support.
Define the operator
$K_{\w}: \mathcal{D}\rightarrow \ell^2(\Dd)$ as
\begin{equation}
K_{\w}(c)=\sum_{i\in I}c_i\T\p_i(\w).
\label{def-operador-k}
\end{equation}
The proof of the following proposition can be found in \cite[Theorem 3.2.3]{christensen}.

\begin{proposition}\label{sintesis-analisis}
The operator $K_{\w}$  defined above is bounded if and only the set
$\{\T\p_i(\w):i\in I\}$ is a Bessel sequence in $\ell^2(\Dd)$.

In that case the adjoint operator of $K_{\w}$, $K^*_{\w}:\ell^2(\Dd)\rightarrow\ell^2(I)$, is given by
$$K^{*}_{\w}(a)=(\langle \T\p_i(\w), a\rangle_{\ell^2(\Dd)})_{i\in I}.$$
\end{proposition}

The operator $K_{\w}$ is called the {\it synthesis operator} and $K^*_{\w}$  the {\it analysis operator}.

\begin{definition}
Let $\{\p_i:i\in I\}\subset L^2(G)$ be a  countable subset and $K_{\w}$  and $K^*_{\w}$ the synthesis and analysis operators. We define the {\it Gramian } of $\{\T\p_i(\w):i\in I\}$ as the operator
$\mathcal{G}_{\w}:\ell^2(I)\rightarrow\ell^2(I)$  given by $\mathcal{G}_{\w}=K^{*}_{\w}K_{\w}$,
and we also define  the {\it dual Gramian } of $\{\T\p_i(\w):i\in I\}$ as the operator
$\mathcal{\tilde{G}}_{\w}:\ell^2(\Dd)\rightarrow\ell^2(\Dd)$  given by $\mathcal{\tilde{G}}_{\w}=K_{\w}K^{*}_{\w}$.
\end{definition}

The Gramian $\mathcal{G}_{\w}$ can be associated with the (possible) infinite matrix
$$\mathcal{G}_{\w}=\Big(\sum_{\dd\in\Dd}\hat{\p_i}(\w+\dd)\overline{\hat{\p_j}(\w+\dd)}\Big)_{i,j\in I}$$
since $\langle \mathcal{G}_{\w}e_i,e_j\rangle=\langle \T\p_i(\w),\T\p_j(\w)\rangle$, where
$\{e_i\}_{i\in I}$ be the standard basis of $\ell^2(I)$.
In a similar way, considering the basis $\{e_{\dd}\}_{\dd\in \Dd}$ of $\ell^2(\Dd)$,
we can associate the dual Gramian $\mathcal{\tilde{G}}_{\w}$  with the matrix
$$\mathcal{\tilde{G}}_{\w}=\Big(\sum_{i\in I}\hat{\p_i}(\w+\dd)\overline{\hat{\p_i}(\w+\dd')}\Big)
_{\dd,\dd'\in \Dd}.$$

\begin{remark}\label{operadores-acotados}
The  operator $K_{\w}$ ($K^{*}_{\w}$) is bounded if and only if
$\mathcal{G}_{\w}$ ($\mathcal{\tilde{G}}_{\w}$) is bounded. In that case we have $\|K_{\w}\|^2=\|K^{*}_{\w}\|^2=\|\mathcal{G}_{\w}\|=\|\mathcal{\tilde{G}}_{\w}\|$.
\end{remark}

Now we will give a characterization of when  $E_H(\A)$ is a frame (Riesz sequence) for $S(\A)$
in terms of the Gramian $\mathcal{G}_{\w}$ and the dual Gramian $\mathcal{\tilde{G}}_{\w}$.

\begin{proposition}
Let $\A=\{\p_i\,:\,i\in I\}\subset L^2(G)$ be a countable set.
Then,
\begin{enumerate}
\item [(1)] The following are equivalent:
\begin{enumerate}
\item[($a_1$)] $E_H(\A)$ is a Bessel sequence with constant $B$.
\item[($b_1$)] $\textrm{supess}_{\w\in\W}\|\mathcal{G}_{\w}\|\leq B$.
\item[($c_1$)] $\textrm{supess}_{\w\in\W}\|\mathcal{\tilde{G}}_{\w}\|\leq B.$
\end{enumerate}
\medskip
\item [(2)] The following are equivalent:
\begin{enumerate}
\item[($a_2$)] $E_H(\A)$ is a frame for $S(\A)$ with constants $A$ and  $B$.
\item[($b_2$)] For almost every $\w\in\W$,
$$A\|a\|^2\leq\langle\mathcal{\tilde{G}}_{\w}a,a\rangle\leq B\|a\|^2,$$
for all $a\in\textrm{span}\{\T\p_i(\w)\,:i\in I\}.$
\item[($c_2$)] For almost every $\w\in\W$,
$$\sigma(\mathcal{\tilde{G}}_{\w})\subset\{0\}\cup[A,B].$$
\end{enumerate}
\item [(3)] The following are equivalent:
\begin{enumerate}
\item[($a_3$)] $E_H(\A)$ is a Riesz sequence for $S(\A)$ with constants $A$ and  $B$.
\item[($b_3$)] For almost every $\w\in\W$,
$$A\|c\|^2\leq\langle \mathcal{G}_{\w}c,c\rangle\leq B\|c\|^2,$$
for all $c\in \ell^2(I)$.
\item[($c_3$)] For almost every $\w\in\W$
$$\sigma(\mathcal{G}_{\w})\subset[A,B].$$
\end{enumerate}
\end{enumerate}
\label{marco-gramiano}
\end{proposition}

\begin{proof}
It follows easily  from Theorem \ref{teo-fibras-marcos}, Theorem \ref{teo-fibras-bases-riesz}, Proposition \ref{sintesis-analisis} and Remark \ref{operadores-acotados}.
\end{proof}

Note that Corollary \ref{coro-ppal-marcos} and Corollary \ref{coro-ppal-riesz} can also be obtained from the previous
proposition.

\begin{definition}
For an $H$-invariant space $V\subset L^2(G)$ we define the  {\it dimension function of  $V$}
as the map $\dim_{V}:\W\rightarrow \N_0$ given by $\dim_{V}(\w)=\dim J(\w)$,
where  $J$ is the range function associated to $V$.
We also define the  {\it spectrum  of $V$} as
$s(V)=\{\w\in\W\,:\, J(\w)\neq 0\}$.
\end{definition}

As in the $\R^d$ case, every  $H$-invariant space can be decomposed into an orthogonal sum of principal $H$-invariant spaces. This can be
easily obtained as a consequence of  Zorn's Lemma as in  \cite{KR}.
The next theorem establishes a decomposition  of $H$-invariant space with additional properties as in \cite{B}.
We do not include its proof since it follows readily from the $\R^d$ case (see \cite[Theorem 3.3]{B}).

\begin{theorem}
Let  us  suppose that $V$ is an  $H$-invariant space of $L^2(G)$.
Then $V$  can be decomposed as an orthogonal sum
$$V=\bigoplus_{n\in\N}S(\p_n),$$
where  $E_H(\p_n)$ is a  Parseval frame for $S(\p_n)$ and
$s(S(\p_{n+1}))\subset s(S(\p_n))$ for all $n\in\N$.

Moreover, $\dim_{S(\p_n)}(\w)=\|\T\p_n(\w)\|$ for all $n\in\N$, and
$$\dim_{V}(\w)=\sum_{n\in\N}\|\T\p_n(\w)\|, \qquad\textrm{a.e.}\,\,\w\in\W.$$
\end{theorem}

\section*{Acknowledgement}

We are grateful to the referees for their useful comments which have led to a significant improvement to the presentation of the article.

\end{document}